\documentclass[10pt]{amsart}
\usepackage{graphicx}
\usepackage{amsmath}
\usepackage{amsthm}
\usepackage{amsfonts}
\usepackage{amssymb}
\usepackage{cite}

\usepackage{epic}
\usepackage{pstricks}

\theoremstyle{plain}
\newtheorem{theorem}{Theorem}

\newtheorem{lemma}[theorem]{Lemma}

\newcommand{\De}{\Delta}

\newcommand{\si}{\sigma}
\newcommand{\de}{\delta}

\renewcommand{\le}{\leqslant}
\renewcommand{\ge}{\geqslant}

\newcommand{\xx}{\mathbf{x}}

\newcommand{\R}{\mathbb{R}}

\newcommand{\N}{\mathbb{N}}

\renewcommand{\P}{\mathcal{P}}
\newcommand{\Q}{\mathcal{Q}}
\newcommand{\RR}{\mathcal{R}}

\newcommand{\card}{\operatorname{card}}

\newcommand{\conv}{\operatorname{conv}}

\newcommand{\sign}{\operatorname{sign}}

\begin{document}


\title[On the Erd\"os-Szekeres theorem]
{An order-refined and generalized version of the Erd\"os-Szekeres theorem on convex polygons} 

\vspace*{-.99cm}

\author{Iosif Pinelis}  
                   
\address{Department of Mathematical Sciences,
Michigan Technological University,
Hough\-ton, MI 49931, 
USA}
\email{ipinelis@mtu.edu}

\date{\today; file convex-poly/szekeres/\jobname}





\maketitle

\setcounter{section}{-1}



A {\em polygon} or, more specifically, $n$-gon, is defined in this paper as any finite sequence $\P=(V_0,\dots,V_{n-1})$ of points 
on the Euclidean plane $\R^2$; the same definition was used in \cite{rob1,rob2,jgeom,angles,elimin,test} and, essentially, \cite{moret}. 
The points $V_0,\dots,V_{n-1}$ are called the {\em vertices} of $\P$.
The 
closed intervals
$[V_i,V_{i+1}]:=\conv\{V_i,V_{i+1}\}$ for $i\in\{0,\dots,n-1\}$
are called here the {\em edges} of polygon $\P$, where 
$V_n:=V_0$.
The symbol $\conv$ denotes, as usual, the convex hull of a set, that is, the intersection of all convex sets containing the given set.
Let us define the convex hull and dimension 
of a polygon $\P$ as, respectively, 
the convex hull and dimension 
of the set of its vertices: $\conv\P:=\conv\{V_0,\dots,V_{n-1}\}$ and $\dim\P:=\dim\{V_0,\dots,V_{n-1}\}=\dim\conv\P$.


Next, define a {\em convex polygon} as a polygon $\P$ such that the union of the edges of $\P$ coincides with the boundary of $\conv\P$. 
Let us emphasize that a polygon in this paper is a sequence and therefore ordered. In particular, even if the vertices $V_0,\dots,V_{n-1}$ of a polygon $\P=(V_0,\dots,V_{n-1})$ are the extreme points of the convex hull of $\P$, it does not necessarily follow that $\P$ is convex. For example, consider the points $V_0=(0,0)$, $V_1=(1,0)$, $V_2=(1,1)$, and $V_3=(0,1)$. Then polygon $(V_0,V_1,V_2,V_3)$ is convex, while polygon $(V_0,V_2,V_1,V_3)$ is not.


Let us say that a polygon $\P=(V_0,\dots,V_{n-1})$ is {\em strict} if 
for any three distinct $i$, $j$, and $k$ in the set $\{0,\dots,n-1\}$, the vertices $V_i$, $V_j$, and $V_k$ are non-collinear. 
Let us say that a polygon is {\em strictly convex} if 
it is both strict and convex. 
Let us say that a polygon $\P=(V_0,\dots,V_{n-1})$ is {\em ordinary} if its vertices are all distinct from one another: for any $i$ and $j$ in the set $\{0,\dots,n-1\}$ such that $i\ne j$, one has   
$V_i\ne V_j$.


For a polygon $\P=(V_0,\dots,V_{n-1})$, let $x_i$ and $y_i$ denote the coordinates of its vertices $V_i$, so that
$V_i=(x_i,y_i)$ for $i\in\{0,\dots,n-1\}$. 
Thus, there is a one-to-one correspondence between all $n$-gons $\P$ and all image points
$$\xx(\P):=(x_1,y_1,\dots,x_{n-1},y_{n-1})$$
in $\R^{2n}$. 
Introduce also the determinants
\begin{equation*}\label{eq:De}
\De_{i,j,k}:=\De_{i,j,k}(\P):=
\left|
\begin{matrix}
1&x_i&y_i \\
1&x_j&y_j \\
1&x_k&y_k 
\end{matrix}\right|
\end{equation*}
for all $i$, $j$, and $k$ in the set $\{0,\dots,n-1\}$. 

In \cite{test}, we proved the following polygon convexity test, which is minimal in a certain exact sense.
\begin{theorem}[\!\!\cite{test}]\label{th:calculation}
A strict $n$-gon $\P=(V_0,\dots,V_{n-1})$ with $n\ge4$ is convex if and only if 
$\sign\De_{i-1,i,i+1}=\sign\De_{0,j,j+1}=\sign\De_{0,1,k+1}$ for all $i$ and $k$ in the set $\{2,\dots,n-2\}$ and all $j$ in
$\{1,\dots,n-2\}$.
\end{theorem}

If $\P=(V_0,\dots,V_{n-1})$ is a polygon, let us refer to any subsequence $(V_{i_0},\dots,V_{i_{k-1}})$ of $\P$, with $0\le i_0<\dots<i_{k-1}\le n-1$,  as a {\em sub-polygon} or, more specifically, 
as a {\em sub-$k$-gon} of $\P$.

In \cite[Corollary~1.24]{elimin}, we proved 
\begin{theorem}[\!\!\cite{elimin}]\label{cor:4-upwards} 
If all sub-$4$-gons of a polygon $\P$ are convex, then $\P$ is convex. 
\end{theorem}

Let us also recall Ramsey's theorem. Let $p$, $q$, and $r$ be any natural numbers such that $p\ge r$ and $q\ge r$, and let $S$ be any finite set. Suppose that every $r$-subset (that is, every subset of cardinality $r$) of the set $S$ is classified arbitrarily as ``good'' or ``bad''. Let us say that a subset $T$ of $S$ is ``$r$-totally good'' if all $r$-subsets of $T$ are good; similarly defined are ``$r$-totally bad'' subsets of $S$. 

\begin{theorem}[Ramsey \cite{ramsey}]\label{th:ramsey} 
There exists either an $r$-totally good $p$-subset of $S$ or an $r$-totally bad $q$-subset of $S$ -- provided that the set $S$ is large enough, that is, provided that $\card S\ge\rho(p,q;r)$, where $\rho(p,q;r)$ is a natural number which depends only on $p$, $q$, and $r$ (but not on the good-bad classification).  
\end{theorem}

The smallest possible lower bound $\rho(p,q;r)$ in Ramsey's theorem is called the Ramsey number and denoted by $R(p,q;r)$.

The main result of this note is
\begin{theorem}\label{th:main}
There is a function $F\colon\N\to\N$ such that for any $k\in\N$ and any 
$n$-gon $\P$ with $n\ge F(k)$ there is a convex sub-$k$-gon of $\P$. 
\end{theorem}

The proof of Theorem~\ref{th:main} is based on 
Theorems~\ref{th:calculation}, \ref{cor:4-upwards}, and \ref{th:ramsey}, stated above, as well as on some other results in \cite{elimin}. The following lemma can be considered a special case of Theorem~\ref{th:main} -- when $k=4$ and polygon $\P$ is strict.

\begin{lemma}\label{lem:13}
Any strict $n$-gon with $n\ge13$ has a convex sub-$4$-gon. 
\end{lemma}

\begin{proof}[Proof of Lemma~\ref{lem:13}]
Let $\P=(V_0,\dots,V_{n-1})$ be any strict $n$-gon with $n\ge13$. 
Let us classify a subset $\{i,j,k\}$ of the set $\{0,\dots,n-1\}$ with $i<j<k$ as ``good'' or ``bad'' according to whether the determinant $\De_{i,j,k}$ 
is positive or negative, respectively.
By Theorem~\ref{th:calculation}, a sub-$4$-gon $(V_{i_0},V_{i_1},V_{i_2},V_{i_3})$ of $\P$ is convex if and only if the index set $\{i_0,i_1,i_2,i_3\}$ is either $3$-totally good or $3$-totally bad. It remains to refer to the result of \cite{13}, which states that the Ramsey number $R(4,4;3)$ is $13$.
\end{proof}

Now one is ready to prove Theorem~\ref{th:main} for any $k$ in the case when the given polygon $\P$ is strict:

\begin{lemma}\label{lem:strict}
There is a function $F\colon\N\to\N$ such that for any $k\in\N$ and any \emph{strict} 
$n$-gon $\P$ with $n\ge F(k)$ there is a convex sub-$k$-gon of $\P$. 
\end{lemma}

\begin{proof}[Proof of Lemma~\ref{lem:strict}]
Let $F(k):=k$ for $k\in\{1,2,3\}$ and $F(k):=R(k,13;4)$ for $k\in\{4,5,\dots\}$. Let $\P:=(V_0,\dots,V_{n-1})$ be a strict $n$-gon with $n\ge F(k)$.
The case $k\in\{1,2,3\}$ is trivial, since all $k$-gons with $k\in\{1,2,3\}$ are convex.
It remains to consider the case $k\in\{4,5,\dots\}$.
Let us classify an $m$-subset $\{i_0,\dots,i_{m-1}\}$ of the set $\{0,\dots,n-1\}$ with $i_0<\dots<i_{m-1}$ 
as ``good'' or ``bad'' according to whether the corresponding sub-$m$-gon $(V_{i_0},\dots,V_{i_{m-1}})$ of $\P$ is convex or not convex, respectively.
Then, by Theorem~\ref{th:ramsey}, there exists either a $4$-totally good
$k$-subset of the set $\{0,\dots,n-1\}$ or a $4$-totally bad $13$-subset of $\{0,\dots,n-1\}$.
However, the second, ``totally bad'' possibility is excluded by Lemma~\ref{lem:13}. 
Hence, all sub-$4$-gons of some sub-$k$-gon of the polygon $\P$ are convex. By Theorem~\ref{cor:4-upwards}, 
such a sub-$k$-gon must be convex. 
\end{proof}

\begin{proof}[Proof of Theorem~\ref{th:main}]
The set (say N\!S\!P($n$)) of all $n$-gons $\P$ that are not strict is the union of all sets of the form $\{\P\colon\De_{i,j,k}(\P)=0\}$ over all triples of integers $(i,j,k)$ such that $0\le i<j<k\le n-1$; hence, the set of image points $\{\xx(\P)\colon\P\in\text{N\!S\!P}(n)\}$ 
is nowhere dense in $\R^{2n}$. Therefore, there exists an infinite sequence $(\P^{(m)})$ of strict $n$-gons $\P^{(m)}=(V_0^{(m)},\dots,V_{n-1}^{(m)})$ such that $\xx(\P^{(m)})\to\xx(\P)$ as $m\to\infty$. 
By Lemma~\ref{lem:strict}, for each $m$ there exist integers $i_0^{(m)},\dots,i_3^{(m)}$ such that $0\le i_0^{(m)}<\dots<i_3^{(m)}\le n-1$ and the sub-$4$-gon $\Q^{(m)}:=(V_{i_0^{(m)}}^{(m)},\dots,V_{i_3^{(m)}}^{(m)})$ of $n$-gon $\P^{(m)}$ is convex. Passing (if necessary) to a subsequence of sequence $(\P^{(m)})$, one has without loss of generality (w.l.o.g.) that the integers $i_0^{(m)},\dots,i_3^{(m)}$ do not depend on $m$, so that $\Q^{(m)}=(V_{i_0}^{(m)},\dots,V_{i_3}^{(m)})$ for some integers $i_0,\dots,i_3$ such that $0\le i_0<\dots<i_3\le n-1$ and all $m$; moreover, 
in view of Theorem~\ref{th:calculation}, the sign of $\de_j(\P^{(m)})$ is w.l.o.g.\ the same 
for all $m$ and all $j\in\{0,1,2,3\}$, where $\de_j(\RR):=\De_{i_0,\dots,i_{j-1},i_{j+1},\dots,i_3}(\RR)$ for any $n$-gon $\RR$. 

It remains to notice that the sub-$4$-gon $\Q:=(V_{i_0},\dots,V_{i_3})$ of $n$-gon $\P$ is convex. Indeed, at least one of the following three cases takes place.

\noindent\emph{Case 1: $\dim\Q\le1$}.\quad Then $\Q$ is convex by \cite[Proposition~1.2]{elimin}.

\noindent\emph{Case 2: $\dim\Q=2$ and $\Q$ is not ordinary}.\quad Then w.l.o.g.\ one has at least one of the following two subcases.

\emph{Subcase 2.1: $V_{i_0}=V_{i_1}$}.\quad Then again $\Q$ is convex, since the convex hull and union of the edges of $4$-gon $\Q=(V_{i_0},\dots,V_{i_3})$ are the same as those of the $3$-gon $(V_{i_1},\dots,V_{i_3})$, and all $3$-gons are convex. 

\emph{Subcase 2.2: $V_{i_0}=V_{i_2}$}.\quad Then 
$$\de_0(\P)=\De_{i_1,i_2,i_3}(\P)=\De_{i_1,i_0,i_3}(\P)=-\De_{i_0,i_1,i_3}(\P)=-\de_2(\P)\ne0$$
(the inequality here is due to the assumptions $\dim\P=2$ and $V_{i_0}=V_{i_2}$).
But this is a contradiction, since $\de_j(\P^{(m)})\to\de_j(\P)$ and the sign of $\de_j(\P^{(m)})$ is the same 
for all $m$ and all $j\in\{0,1,2,3\}$.

\noindent\emph{Case 3: $\Q$ is ordinary}.\quad According to \cite[Proposition~1.9]{elimin}, for each $m$ the convex polygon $\Q^{(m)}$ is to-one-side. By Definition~\cite[Definition~1.8]{elimin}), this means that for each $m$ and each $j\in\{0,1,2,3\}$ there exists a linear functional $\ell_{m,j}$ of norm $1$ on the Euclidian space $\R^2$ such that $\ell_{m,j}(V_{i_{j+1}}^{(m)})=\ell_{m,j}(V_{i_j}^{(m)})$, while $\ell_{m,j}(V_{i_r}^{(m)})\ge\ell_{m,j}(V_{i_j}^{(m)})$ for all $r\in\{0,1,2,3\}$ \big(here for $j=3$ the point $V_{i_{j+1}}^{(m)}$ is to be understood as $V_{i_0}^{(m)}$\big). By the compactness of the unit circle, w.l.o.g.\ the limits $\ell_j:=\lim_{m\to\infty}\ell_{m,j}$ exist and are of norm $1$, for all $j\in\{0,1,2,3\}$. Since $V_i^{(m)}\to V_i$ for all $i$, one concludes that the ``limit'' $4$-gon $\Q=(V_{i_0},\dots,V_{i_3})$ is to-one-side. Using \cite[Proposition~1.9]{elimin} again, one sees that $\Q$ is quasi-convex. Finally, the condition that $\Q$ is ordinary and \cite[Proposition~1.13]{elimin} imply that $\Q$ is convex.

\end{proof}

If $\si\colon\{0,\dots,n-1\}\to\{0,\dots,n-1\}$ is a bijection, let us refer to the polygon $\si(\P):=(V_{\si(0)},\dots,V_{\si(n-1)})$ as a \emph{permutation of polygon} $\P=(V_0,\dots,V_{n-1})$. Let us say that a polygon $\P$ is \emph{pre-convex} if a permutation of $\P$ is convex. By \cite[Proposition~1.5]{elimin}, a strict polygon $\P$ is pre-convex iff the set of all vertices of $\P$ coincides with the set of all extreme points of $\conv\P$. The polygon convexity is a much more restrictive notion than that of the polygon pre-convexity. In particular, in view of the main result of \cite{angles}, any cyclic polygon is pre-convex. On the other hand, by \cite[Proposition~1.6]{elimin}, only $2n$ permutations of the $n!$ permutations of any given strictly convex $n$-gon $\P$ with $n\ge3$ are convex (while all the $n!$ permutations of $\P$ are obviously pre-convex). 

Theorem~\ref{th:main} with the term ``convex'' replaced by ``pre-convex'' and under the additional assumption that the given polygon is strict is the famous Erd\"os-Szekeres theorem \cite{szek}, which has been followed by a great many generalizations; see e.g.\ survey \cite{survey}. \big(In \cite{szek} and in many other papers, the term ``convex'' was used essentially in place of ``pre-convex'' -- but without the notion of polygon convexity or that of a polygon itself being formally defined.\big) Thus, Theorem~\ref{th:main} refines and generalizes the Erd\"os-Szekeres theorem. In particular, this addresses the comment made in \cite[page~464]{szek} that ``It is desirable to extend the usual definition of convex polygon to include the cases where three or more consecutive points lie on a straight line.''

An open problem that remains is to determine, for each $k$, the least possible number $F(k)$ in Theorem~\ref{th:main}. In particular, one can ask whether $13$ in Lemma~\ref{lem:13} can be replaced by a smaller number; that is, whether $F(4)<13$. 
The example of the $7$-gon
$\big((-13,0),(15,0),(0,16),(18,39),(27,-15),(10,20),(16,30)\big)$, none of whose $\binom 74=35$ sub-$4$-gons is convex, shows that $F(4)\ge8$.   

\renewcommand{\refname}{\textsf{\bf Literature}}

\end{document}